 \documentclass[12pt]{amsart}
 \usepackage{amsfonts}

 \usepackage{epsfig, graphicx, graphics, color, amsmath,amssymb, amsfonts, amstext,amsthm, latexsym}

   \newtheorem{lemma}{Lemma}[section]
   \newtheorem{theorem}[lemma]{Theorem}
   \newtheorem{remark}[lemma]{Remark}

   \newtheorem{definition}[lemma]{Definition}

   \newcommand{\EX}{{\Bbb{E}}}
   \newcommand{\PX}{{\Bbb{P}}}
\newcommand{\s}{\sigma}

\newcommand{\Om}{\Omega}

\newcommand{\p}{\partial}

\renewcommand{\phi}{\varphi}
\newcommand{\N}{{\mathbb N}}
\newcommand{\R}{{\mathbb R}}

\newcommand{\F}{{\mathcal{F}}}
\newcommand{\B}{{\mathcal{B}}}

\makeatletter
    \newcommand\figcaption{\def\@captype{figure}\caption}
    \newcommand\tabcaption{\def\@captype{table}\caption}
\makeatother

\title[Sample Spaces for Random Dynamical Systems]
{Canonical Sample Spaces\\ for Random Dynamical Systems}

\author{Jinqiao Duan }
\address[Jinqiao Duan]
{Department of Applied Mathematics\\
 Illinois Institute of Technology\\
   Chicago, IL 60616, USA }
\email[J.~Duan]{duan@iit.edu}

\author{Xingye Kan}
\address[Xingye Kan]
{Department of Applied Mathematics\\
 Illinois Institute of Technology\\
   Chicago, IL 60616, USA }
 \email[X. Kan]{xkan@iit.edu}

\author{Bj{\"o}rn Schmalfu{\ss }}
\address[Bj{\"o}rn Schmalfu{\ss }]{Institut f\"{u}r Mathematik\\
Fakult\"{a}t EIM, Warburger Stra{ss}e 100, 33098\\
Paderborn, Germany
 }
\email[Bj{\"o}rn Schmalfu{\ss }]{schmalfuss@upb.de}

 \date{\today}

\keywords{Random dynamical systems, Brownian motion, fractional brownian motion,
Wiener measure, L\'evy motion, colored noise. \\
This work was partially supported by NSF grants 0620539 and
0731201, and a Graduate Research Fellowship from Illinois
Institute of Technology}

\begin{document}
\def \R {{\bf R}}
\def \C {{\bf C}}
\def \E {{\bf E}}
\def \N {{\bf N}}
\def \P {{\bf P}}
\def \Q {{\bf Q}}
\def \l{\langle}
\def \r{\rangle}
\def \p{\partial}
\def \ba{\begin{eqnarray*}}
\def \ea{\end{eqnarray*}}
\def \bc{\begin{center}}
\def \ec{\end{center}}
\def \bl{\begin{flushleft}}
\def \el{\end{flushleft}}
\def \br{\begin{flushright}}
\def \er{\end{flushright}}
\def \bd{\begin{document}}
\def \ed{\end{document}}

\begin{abstract}
This is an overview about   natural  sample spaces for
differential equations driven by various noises. Appropriate
sample spaces are needed in order to facilitate a random dynamical
systems approach for stochastic differential equations. The noise
could be white or colored, Gaussian or non-Gaussian, Markov or
non-Markov, and semimartingale or non-semimartingale.  Typical
noises are defined in terms of Brownian motion, L\'evy motion and
fractional Brownian motion. In each of these cases, a canonical
sample space with an appropriate metric (or topology that gives
convergence concept) is introduced. Basic properties of
canonical sample spaces, such as separability and completeness,
are then discussed.

Moreover, a flow defined by shifts, is introduced on these
canonical sample spaces. This flow has an invariant measure which
is the probability distribution for Brownian motion, or L\'evy
motion or fractional Brownian motion.
Thus canonical sample spaces are much richer in mathematical
structures than the usual sample spaces in probability theory, as
they have metric or topological structures, together with a
  shift flow (or driving flow) defined on it.  This
facilitates dynamical systems approaches for studying stochastic
differential equations.

\end{abstract}

\maketitle

\section{Random dynamical systems}
\label{RDS}

Stochastic differential equations (SDEs) or stochastic partial
differential equations (SPDEs) arise as mathematical models for
complex systems under various random influences in engineering and
science. Here we only consider random dynamical systems defined by
SDEs. Such SDEs define random dynamical systems (RDS) with
appropriate sample spaces, much as ordinary differential equations
define deterministic dynamical systems.

Theory of random dynamical systems allows to discuss the
qualitative behavior of stochastic systems that are not only
driven by a white noise, Markov processes and semimartingales, but
also driven by non-Markov processes or by non-semimartingales
(e.g., fractional Brownian motion). To analyze these more general
noise cases, appropriate sample spaces and  ergodic theory play an
important role. In this article, we discuss canonical or natural
sample spaces for SDEs with various noises.

We recall the definition of a random dynamical system (RDS) in the
state space $H=\R^n$, with the underlying  probability space
$(\Om, \F, \PX)$, and   with time $t$ varying in
$\mathbb{T}=\R=(-\infty, \infty)$ or $\mathbb{T}=\R^+=[0,
\infty)$, as in Arnold \cite{Arnold}. The state space $\R^n$ is
equipped with the Euclidean norm (or length) $|x|
=\sqrt{x_1^2+\cdots +x_n^2}$ and the usual scalar product
$<x,y>=x_1y_1+\cdots +x_ny_n$.

Note that a deterministic dynamical system on the state space $H$
is a mapping $\psi: \mathbb{T}\times H \to H, \; (t, x) \mapsto
\psi(t,x)$, such that   the flow property is satisfied:
$$
\psi(0,x)=x,\;\; \psi(t+s, x)=\psi(t, \psi(s,x)),
$$
for all $t,s \in \mathbb{T} $ and $x \in H$.

For a random dynamical system, we need an extra ingredient,
namely, a model for the noise. Moreover, the flow property has a
twist (thus called cocycle property) due to the effect of noise.

 \begin{definition} (Random dynamical system)\\
 A   random dynamical system (RDS), denoted by
$\varphi$, consists of two ingredients: \\
(i) \textbf{Model for the noise}:   A driving   flow on a
probability space $(\Om, \F, \PX)$, i.e., a flow $(\theta_t)_{t\in
\mathbb{T}}$ on the sample space $\Omega$, such that $\mathbb{P}$
is invariant, namely $\theta_t \mathbb{P} = \mathbb{P}$ for all $t
\in \mathbb{T}$, and $(t,\omega)\mapsto
\theta_t \omega$ is measurable from $\mathbb{T} \times \Om$ to $\Om$. \\
(ii) \textbf{Model for the evolution}:  A cocycle $\varphi$ over
$\theta$, i.e. a measurable mapping $\varphi: \mathbb{T}
\times\Omega\times H\rightarrow  H, $ $(t,\omega,x)\mapsto
\varphi(t,\omega,x)$, such that   the family
$\varphi(t,\omega,\cdot) = \varphi(t,\omega): H \rightarrow H$ of
random  mappings  satisfies the cocycle property:
\begin{eqnarray}
\varphi(0,\omega)=id_H , \varphi(t+s,\omega)=\varphi(t,\theta_s
\omega)\circ \varphi(s,\omega) \mbox{ for all } t,s \in
\mathbb{T},
 \omega\in \Omega. \label{5}
\end{eqnarray}
\end{definition}
Here the driving flow $\theta_t$ describes stationary dynamics of
noise in an appropriately chosen sample space (see below). The
mathematical model for noises in engineering and science is
usually a stationary generalized stochastic process \cite{Hida74}.

 When  $(t,x) \mapsto \varphi(t,\omega,x)$ is continuous for
all $\omega \in \Omega$, we say that $\varphi$ is a continuous
RDS. Since   the continuity in space $x$ is quite common for RDS
generated by stochastic differential equations, we usually do not
specifically mention this spatial continuity. We often call
$\varphi$ a continuous-time or discrete-time RDS when it is
continuous or discrete in time $t$.

It follows from  \cite{Arnold}  that $\varphi(t,\omega)$, $t\in
\R$,  is a homeomorphism of $H$ and
 $$\varphi(t,\omega)^{ -1 }=\varphi(-t,\theta_t\omega).$$

\section{Dynamical systems driven by white noises}
\label{whitenoise}

\subsection{Brownian Motion}

The physical phenomenon \emph{Brownian motion}\footnote{The
phenomenon was first observed by Jan Ingenhouz in 1785, but was
subsequently rediscovered by Brown in 1828, according to sources
used by Eric Weisstein's \emph{World of Physics}, which can be
found on the Internet at
http://scienceworld.wolfram.com/physics/BrownianMotion.html} is
due to the incessant hitting of pollen by the much smaller
molecules of the liquid. The hits occur a large number of times in
any small time interval, independently of each other and the
effect of a particular hit is small compared to the total effect
\cite{Plato}. The physical theory of this motion, set up by Albert
Einstein in 1905, suggests
that the motion is random, and has the following properties:\\
\\
i) the motion is continuous;\\
ii) it has independent increments;\\
iii) the increments are stationary and Gaussian random variables.\\

Figure 1   shows a sample path of the Brownian motion.
\begin{figure} \label{BMpath}
 \includegraphics[height=3in,width=4in]{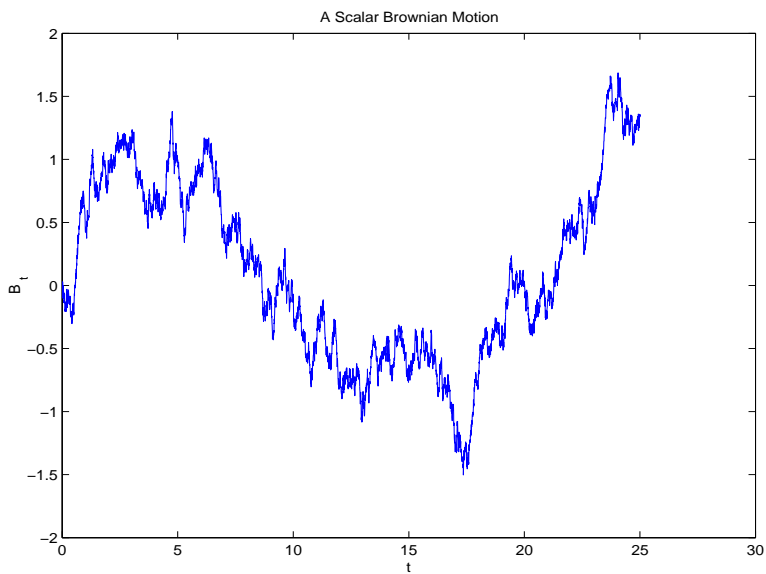}
\caption{\small \sl A sample path of   Brownian motion $B(t)$}
\end{figure}

Intuitively speaking, property i) says that the sample path of the
Brownian motion is continuous. Property ii) means that the
displacements of a pollen particle over disjoint time intervals are
independent random variables. Property iii) is natural considering
the \emph{Central Limit Theorem}.\\
\\
We now describe the Brownian motion in the mathematical language,
i.e., introduce the first definition of the Brownian
motion\cite{Hida74, Arnold74, BrzezniakZastawniak}.\\
\\
\noindent\textbf{Definition A:} A stochastic process
$\{B_{t}(\omega):t\geq0\}$ defined on a probability space
$(\Omega,\mathcal {F},P)$ is called a \emph{Brownian motion} or a
\emph{Wiener process} if   the following conditions
hold:\\
1) $B_{0}(\omega)=0$ a.s.; \\
2) the sample paths $t\rightarrow B_{t}(\omega)$ are a.s. continuous;\\
3) $B_{t}(\omega)$ has stationary independent increments;\\
4) the increments $B_{t}(\omega)-B_{s}(\omega)$ has the normal
distribution with mean $0$ and variance $t-s$, i.e.
$B_{t}(\omega)-B_{s}(\omega)\sim N(0,t-s)$ for any $0\leq s <
t$.\\
\\
Since the stochastic process $B_{t}(\omega)$ is a mapping from a
probability space to a metric space and is governed by its law, i.e.
the probability measure on the metric space induced by
$B_{t}(\omega)$, we want to find it in order to have a better
understanding of the Brownian motion. In fact, we can find the
finite dimensional distribution of the Brownian motion using
condition 3) and 4) in Definition A, and this gives another
definition of the Brownian motion as we shall see. We first write
down the second definition\cite{BrzezniakZastawniak} and then prove it is equivalent to Definition A.\\
\\
\noindent\textbf{Definition B:} A stochastic process
$\{B_{t}(\omega):t\geq0\}$ defined on a probability space
$(\Omega,\mathcal {F},P)$ is called a \emph{Brownian motion} or a
\emph{Wiener process} if   the following conditions
hold:\\
1') $B_{0}(\omega)=0$ a.s.; \\
2') the sample paths $t\rightarrow B_{t}(\omega)$ are a.s. continuous;\\
3') for any finite sequence of times
$0=t_{0}<t_{1}<t_{2}<\cdots<t_{n}$ and Borel sets
$B_{1},\cdots,B_{n}\subset\mathbb{R}$
\begin{align}\label{eq1}
P\{B_{t_{1}}(\omega)\in B_{1},\cdots,B_{t_{n}}(\omega)\in B_{n}\} \nonumber \\
=\int_{B_{1}}\cdots
\int_{B_{n}}p(t_{1},0,x_{1})p(t_{2}-t_{1},x_{1},x_{2})\cdots\nonumber \\
\cdots p(t_{n}-t_{n-1},x_{n-1},x_{n})dx_{1}\cdots dx_{n}
\end{align}
where
\begin{equation}
p(t;x,y)=\frac{1}{\sqrt{2\pi t}}e^{-\frac{(x-y)^{2}}{2t}}
\end{equation}
defined for any $x,y\in\mathbb{R}$ and $t>0$ is called the
\emph{transition density}.\\
\\
We can see the conditions 1'),2') in Definition B are completely the
same as 1),2) in Definition A, respectively, the only thing we need
to do in the proof of the equivalence of the two is to show 3),4)
$\Leftrightarrow$ 3').\\
\\
\noindent\textbf{Proof} 3),4)$\Rightarrow$ 3')\\
Firstly, we show that Brownian motion $B_{t}$ has Markov property,
by proving that the conditional distribution of $B_{t+s}$ given
$\mathcal{F}_t$ is the same as that given $B_t$, in terms of
moment generating function \cite{Klebaner}. In fact,\\
\begin{align*}
E(e^{u B_{t+s}}|\mathcal {F}_{t})
&=E(e^{u[(B_{t+s}-B_{t})+B_{t}]}|\mathcal {F}_{t}) \\
&=e^{uB_{t}}E(e^{u(B_{t+s}-B_{t}}|\mathcal {F}_{t})\\
&=e^{uB_{t}}E(e^{u(B_{t+s}-B_{t})})\\
&=e^{uB_{t}}e^{u^{2}s/2}\\
&=e^{uB_{t}}E(e^{u(B_{t+s}-B_{t})}|B_{t})\\
&=E(e^{uB_{t+s}}|B_{t}).
\end{align*}
\\
Secondly, we compute the joint distribution of the Brownian motion.
\begin{align*}
&P\{B(t_{k+1})\leq x_{k+1},B(t_{k})\leq x_{k}\}\\
&=P\{B(t_{k})\leq x_{k},[B(t_{k+1})-B(t_{k})]+B(t_{k})\leq x_{k+1}\}\\
&=P\{B(t_{k})\leq x_{k},[B(t_{k+1})-B(t_{k})]\leq x_{k+1}-B(t_{k})\}\\
&=P\{B(t_{k+1})-B(t_{k})\leq x_{k+1}-B(t_{k})|B(t_k)\leq x_k\}P\{B(t_{k})\leq x_{k}\}\text{ (Conditional probability)}\\
&=\int_{-\infty}^{x_{k}}\int_{-\infty}^{x_{k+1}-x}\frac{1}{\sqrt{2\pi(t_{k+1}-t_{k})}}e^{-{\frac{y^{2}}{2(t_{k+1}-t_{k})}}}dy
\frac{1}{\sqrt{2\pi t_{k}}}e^{-{\frac{x^{2}}{2t_{k}}}}dx\\
&=\int_{-\infty}^{x_{k}}\frac{1}{\sqrt{2\pi
t_{k}}}e^{-{\frac{x^{2}}{2t_{k}}}}dx\int_{-\infty}^{x_{k+1}}\frac{1}{\sqrt{2\pi(t_{k+1}-t_{k})}}e^{-{\frac{(y-x)^{2}}{2(t_{k+1}-t_{k})}}}dy\\
&=\int_{-\infty}^{x_{k}}p(t_{k};0,x)dx\int_{-\infty}^{x_{k+1}}p(t_{k+1}-t_{k};x,y)dy.
\end{align*}
Finally, we have the finite dimensional distribution.
\begin{align*}
&\qquad P\{B(t_{1})\leq x_{1},B(t_{2})\leq x_{2},\cdots,
B(t_{n})\leq
x_{n}\}\\
&\ =P\{B(t_{1})\leq x_{1}|B(t_{0})=0\}P\{B(t_{2})\leq
x_{2}|B(t_{1})\leq x_{1}\}\\
&\qquad P\{B(t_{3})\leq x_{3}|B(t_{1})\leq x_{1},B(t_{2})\leq
x_{2}\} \cdots P\{B(t_{n})\leq x_{n}|B(t_{1})\leq x_{1}\cdots
B(t_{n-1})\leq x_{n-1}\}\\ &\ =P\{B(t_{1})\leq
x_{1}|B(t_{0})=0\}P\{B(t_{2})\leq
x_{2}|B(t_{1})\leq x_{1}\}\\
&\qquad P\{B(t_{3})\leq x_{3}|B(t_{2})\leq x_{2}\} \cdots
P\{B(t_{n})\leq
x_{n}|B(t_{n-1})\leq x_{n-1}\}\text{ (Markov property)}\\
&\ =\frac{P\{B(t_{1})\leq x_{1}\}P\{B(t_{2})\leq x_{2},B(t_{1})\leq
x_{1}\}\cdots P\{B(t_{n})\leq x_{n},B(t_{n-1})\leq
x_{n-1}\}}{P\{B(t_{1})\leq x_{1}\}P\{B(t_{2})\leq x_{2}\}\cdots
P\{B(t_{n-1})\leq x_{n-1}\}}\\
&\
=\int_{-\infty}^{x_{1}}p(t_{1};0,y_{1})dy_{1}\int_{-\infty}^{x_{2}}p(t_{2}-t_{1};y_{1},y_{2})dy_{2}\cdots
\int_{-\infty}^{x_{n}}p(t_{n}-t_{n-1};y_{n-1},y_{n})dy_{n}
\end{align*}
which is obviously the same as
3').\qquad\qquad\qquad\qquad\qquad\qquad\qquad\qquad\qquad\qquad\qquad\qquad\qquad\qquad
$\Box$\\
Conversely, we need to show 3')$\Rightarrow$ 3),4)\\
However, this part of work is completely done in
\cite{BrzezniakZastawniak}(see page 153-155).\\

\subsection{\bf Wiener Measure}\label{s2}

In Section 1, we treat the Brownian motion as a stochastic
process, i.e., a collection of time-parameterized random
variables. Since a stochastic process $\xi$ is governed by its
law, i.e., the probability measure $\mu:=P\xi^{-1}$ on the space
it maps to, one may ask the question why the Brownian motion can
be determined only by its finite dimensional distributions
although we have proved the so-defined Bownian motion (Definition
B) coincides with the definition describing the phenomenon
(Definition A). Motivated by this question, we shall treat the
Brownian motion as a ``random variable", which we call
\emph{random function} \cite{Billingsley}, and this point of view
introduces the \emph{Wiener measure}, the probability measure
defined on an appropriate space which gives us the Brownian
motion. We will sketch the ideas showing the existence and
uniqueness of this special probability measure, starting from
three different spaces, i.e, $\mathbb{R}[0,\infty), C[0,\infty),
C_{0}[0,\infty)$, the spaces of arbitrary functions, continuous
functions, and continuous functions passing zero at time $0$, defined on $[0, \infty)$,
respectively.\\
\\
\textbf{First approach} This way of showing the existence and
uniqueness of the Wiener measure is the one most often used to
construct a Markov Process and is rather technical. The main idea
is the following: First   define a set function on the algebra
generated by the cylinder sets in $\mathbb{R}[0,\infty)$ according
to the finite-dimensional distribution of Brownian motion. Note
that it is a set function rather than a probability measure since
it is just finitely additive (not countably additive); then by the
celebrated \emph{Kolmogorov extension theorem} \cite{Hida74}, this
set function is uniquely extended to a probability measure on the
$\sigma-$algebra generated by the algebra mentioned above.\\
\\
More precisely, we now give the definitions and theorems.\\
\\
$\bullet$ $\mathbb{R}[0,\infty)$ denotes the set of all real-valued
functions
on $[0,\infty)$.\\
\\
$\bullet$ \textbf{Cylinder set} \cite{Hida74} A subset of
$\mathbb{R}[0,\infty)$ of the form
\begin{align}\label{eq5}
A=\{\omega\in\mathbb{R}[0,\infty):(\omega(t_{1}),\ldots,\omega(t_{n}))\in
B_{n}\}
\end{align}
where $B_{n}$ is a Borel subset of $\mathbb{R}^{n}$ is called a
\emph{cylinder set}.\\
Fixing $t_{1},\ldots t_{n}$ but varying $B_{n}$ over the entire
Borel subsets of $\mathbb{R}^{n}$, the class of such all cylinder
sets forms a $\sigma-$algebra
$\mathfrak{B}^{(t_{1},\ldots,t_{n})}$.\\
\\
$\bullet$ \textbf{Set function} A set function
$\Phi_{t_{1},\ldots,t_{n}}$ is defined on the measurable space
$(\mathbb{R}[0,\infty),\mathfrak{B}^{(t_{1},\ldots,t_{n})})$, given
by
\begin{align*}
\Phi_{t_{1},\ldots,t_{n}}(A):=P\{(B_{t_{1}},\ldots,B_{t_{n}})\in
B_{n}\}
\end{align*}
By varying the choice of the finite time points
$\{t_{1},\ldots,t_{n}\}\subset[0,\infty)$, we get a class of such
set functions $\Phi:=\{\Phi_{t_{1},\ldots,t_{n}}\}$, and this class
is independent of the cylindrical expression of the functions in
$\mathbb{R}[0,\infty)$, i.e. if the cylinder set $A$ of \eqref{eq5}
has another expression, say
\begin{align*}
A=\{\omega\in\mathbb{R}[0,\infty):(\omega(s_{1}),\ldots,\omega(s_{m}))\in
B_{m}\},
\end{align*}
then
\begin{align*}
\Phi_{t_{1},\ldots,t_{n}}(A)=\Phi_{s_{1},\ldots,s_{m}}(A).
\end{align*}
\\
$\bullet$ \textbf{Algebra and $\sigma$-algebra}\\
$\cdot$ $\mathfrak{U}[0,\infty)$ denotes the algebra of subsets of
$\mathbb{R}[0,\infty)$ consisting all cylinder sets\\
$\cdot$ $\mathfrak{B}[0,\infty)$ denotes the smallest
$\sigma-$algebra containing $\mathfrak{U}[0,\infty)$\\
$\cdot$ $\Phi$ is a finitely additive measure on
$(\mathbb{R}[0,\infty),\mathfrak{U}[0,\infty))$ such that the
restriction of $\Phi$ to $\mathfrak{B}^{(t_{1},\ldots,t_{n})}$
coincides with $\Phi_{t_{1},\ldots,t_{n}}$
\begin{theorem}[\textbf{Kolmogorov extension theorem}]
The set function $\Phi$ on
$(\mathbb{R}[0,\infty),\mathfrak{U}[0,\infty))$ is uniquely
extendable to a probability measure $\tilde{\Phi}$ on
$(\mathbb{R}[0,\infty),\mathfrak{B}[0,\infty))$.
\end{theorem}
It has been proved \cite{Hida74, KaratzasShreve} that there exist
a unique probability measure $P$ on
$(\mathbb{R}[0,\infty),\mathfrak{B}[0,\infty))$, under which the
coordinate mapping process
\begin{align*}
B_{t}(\omega):=\omega(t);\ \omega\in\mathbb{R}[0,\infty),t\geq0
\end{align*}
satisfies condition 3) and 4) of Definition A in section 2.1, so it
does not introduce a ``Brownian motion" as we defined. Fortunately,
by another famous theorem of Kolmogorov\cite{Oksendal}, the
continuity problem has been solved.
\begin{theorem}[\textbf{Kolmogorov continuity theorem}]
Suppose that the process $X=X\{t\}_{t\geq0}$ satisfies the following
condition: for all  $T>0$ there exist positive constants
$\alpha,\beta,D$ such that
\begin{align*}
E[|X_{t}-X_{s}|^{\alpha}]\leq D\cdot|t-s|^{1+\beta};\ 0\leq s,t\leq
T
\end{align*}
then there exists a continuous \emph{modification} of $X$.
\end{theorem}
\textbf{Modification }\cite{Oksendal} Suppose that
$\{X_{t}\},\{Y_{t}\}$ are stochastic processes on $(\Omega,\mathcal
{F},P)$, then we say that $\{X_{t}\}$ is a \emph{modification} of
$\{Y_{t}\}$ if
\begin{align*}
P(\{\omega;X_{t}(\omega)=Y_{t}(\omega)\})=1\ \forall t
\end{align*}
Note that if $X_{t}$ is a modification of $Y_{t}$, then they have
the same finite-dimensional distributions.\\

Now there is one problem needs to be solved: Why is
$B_{0}(\omega)=\omega(0)=0$ a.s.?\\
\\

\noindent\textbf{Definition} \cite{Billingsley} \emph{Wiener
measure}, denoted by $\mu_{W}$, is a probability measure on a
measurable space
$(C,\mathcal {C})$ having the following two properties:\\
i) each $\omega(t)\in C$ is normally distributed under $\mu_{w}$
with mean $0$ and variance $t$, i.e,
\begin{equation*}
\mu_{w}\{\omega(t)\leq x\}=\frac{1}{\sqrt{2\pi
t}}\int_{-\infty}^{x}e^{-{\frac{u^{2}}{2t}}}du, \;  x\in \mathbb{R}.
\end{equation*}
For $t=0$ this is interpreted to mean that
$\mu_{w}\{\omega(0)=0\}=1$.\\
ii) the stochastic process $\{\omega(t): t\geq 0\}$ has independent
increments under $\mu_{w}$, i.e., $\forall 0\leq t_{0}\leq t_{1}\leq
\cdots\leq t_{n}$,
$\omega(t_{1})-\omega(t_{0}),\omega(t_{2})-\omega(t_{1}),\cdots\omega(t_{n})-\omega(t_{n-1})$
are independent under $\mu_{w}$.\\

\textbf{Remark} It may not seem so obvious   that this approach in
fact proves the existence of the Wiener measure. Also, we cannot
obtain the Wiener measure simply by assigning measure one to
$C[0,\infty)$; see \cite{KaratzasShreve}. One might hope to
construct the measure directly on $C[0,\infty)$. Indeed, this is
the main idea of
the second approach we are going to present.\\
\\
\textbf{Second approach} There are some advantages if we start
from $C[0,\infty)$ since we can make it \emph{Polish} (complete
and separable) by assigning an appropriate \emph{metric}. It has
been proved that the Borel $\sigma-$algebra generated by the open
sets in $C[0,\infty)$ is equal to the $\sigma-$algebra generated
by all the cylinder sets \cite{KaratzasShreve}. Thus there is a
totally
different way of constructing Wiener measure.\\
\\
\emph{\textbf{Metric on }$C[0,\infty)$}
\begin{equation}\label{eq2}
    \rho(\omega_{a},\omega_{b}):=\sum_{n=1}^{\infty}\frac{1}{2^{n}}\max_{0\leq t\leq
    n}(|\omega_{a}(t)-\omega_{b}(t)|\wedge 1)
\end{equation}
This metric introduces the topology of uniform convergence on
compact intervals. For convenience, we denote $C[0,\infty)$ by $C$
and its Borel $\sigma$-algebra by $\mathcal {C}$.\\
\\
The main idea is to construct a sequence of probability measure
$\{P_{n}\}$ on $(C,\mathcal {C})$ such that
$P_{n}\Rightarrow\mu_{w}$. Since \eqref{eq1} gives the f.d.d of
the Wiener measure $\mu_{w}$, we must let the f.d.d of $\{P_{n}\}$
weakly converges to those of $\mu_{w}$. Although weak convergence
in $C$ need not follow the weak convergence of the f.d.d alone in
general, it does if we add the condition that $\{P_{n}\}$ is tight
\cite{KaratzasShreve, Billingsley}. In fact, if a metric space is
separable and complete, then each probability measure on the
measurable space is tight (see Theorem 1.3 in \cite{Billingsley}).
It can be shown that equipped with the metric defined by
\eqref{eq2}, the canonical space is separable and complete
\cite{Billingsley}.
For detail, see \cite{KaratzasShreve}\cite{Billingsley}.\\
\\
\textbf{Remark}  What is still not natural is that
$B_{0}(\omega)=\omega(0)=0$ a.s. So next we will talk about
constructing Wiener measure on $C_{0}[0,\infty)$. It appears that
the above two different approaches may be adapted to this case.


We can define the Brownian motion $B_t$ for $t \in \R$ as follows:
Taking two independent Brownian motions $\hat{B}_t$ and
$\tilde{B}_t$, we define
\begin{equation}
  B_t = \begin{cases}
         \hat{B}_t,   &\text{if $t \geq 0$;}\\
         \tilde{B}_{-t},    &\text{if $t<0$.}
    \end{cases}.
\end{equation}
We can   work on the   space $C_0(\R, \R^n)$ for two-sided
Brownian motion. The Wiener measure defined above should be
similar defined in this space \cite{Arnold}.

\subsection{Canonical sample space}

We consider a SDE
\begin{equation}
dX_t   =  b(X_t )d t + \s(X_t )dB(t).
\end{equation}
The canonical sample space is $\Om:=C_0(\R^, \R^n)$, space of
continuous functions that are zero at time zero, equipped with the
compact open topology,     the Borel $\s-$field $\F:=\B(C_0(\R^+,
\R^n))$, and Wiener (probability) measure $\mu_W$.

We   introduce a driving  flow $\theta_t$ on this canonical sample
space $\Om$ is given by the Wiener shifts
 \begin{equation*}
     \theta_t\omega(\cdot)=\omega(\cdot+t)-\omega(t),\quad
     t\in\mathbb{R},\qquad \omega\in\Omega=C_0(\mathbb{R},\mathbb{R}^d).
 \end{equation*}
 In this case the measure $\mu_W$ is invariant \cite{Arnold},
 i.e.,
 \begin{equation*}
\mu_{W}(\theta_{t}^{-1}(A))=\mu_{W}(A)
\end{equation*}
for all $A \in \mathcal{F}$. Moreover,  this invariant measure is
actually
  {\em ergodic} with
 respect to the flow $\theta$ \cite{Boxler}.

\section{Dynamical systems driven by colored noises}
\label{colorednoise}

Colored noise, or noise with non-zero correlation ('memory') in
time, are common in the physical, biological and engineering
sciences \cite{Hanggi3}. A good candidate for modeling colored
noise is the fractional Brownian motion.

\subsection{Fractional Brownian motion }

 A fractional Brownian motion (fBM)
$B^{H}(t)$, $t\in \R$,  with $H\in(0,1)$   the Hurst parameter, is
still a Gaussian process. But it is characterized by the
stationarity of its increments and a  memory property.  The
increments of the fractional Brownian motion are not independent,
except in the standard Brownian motion case ($H=\frac12$). Thus it
is not a  Markov process except when $H=\frac12$. Specifically,
$B^{H}\left( 0\right) =0$ a.s.,   mean $\EX B^{H}\left( t\right)
=0$, covariance   $\EX[B^H(t)B^H(s)] = \frac12(|t|^{2H} + |s|^{2H}
- |t-s|^{2H})$, and variance $Var\left[ B^{H}\left( t\right)
-B^{H}\left( s\right) \right] =\left\vert t-s\right\vert ^{2H}$.
It also exhibits power scaling and path regularity properties with
H\"older parameter $H$, which are very distinct from Brownian
motion. The standard Brownian motion is a special fBM with
$H=1/2$. Figure 2 is a sample path of the fractional Brownian
motion with $H=0.25$.

\begin{figure}[htbp]
\begin{center}
\includegraphics[height=3in,width=4in]{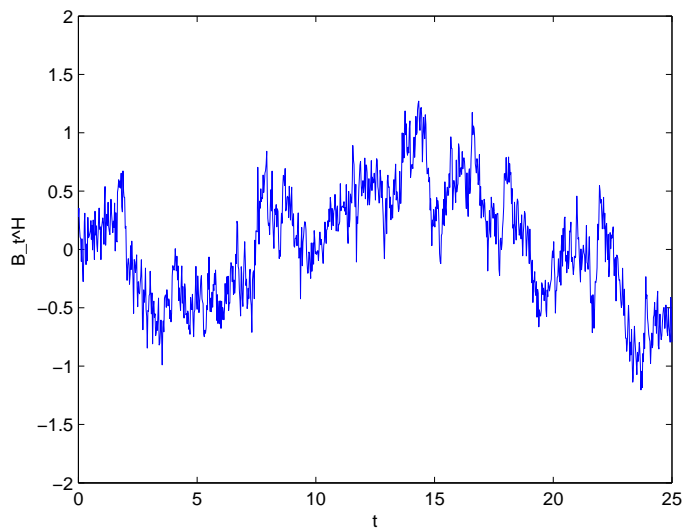}
\caption{\small \sl A sample path of fractional Brownian motion
$B^H(t)$, with $H=0.25$}
\end{center}
\end{figure}

Note that fBM  $B^{H}(t)$ is non-Markov and non-semimartingale.
Thus the usual stochastic integration \cite{Protter} is not
applicable, and other integration concepts are needed
\cite{Nualart, Mishura}.

\subsection{Canonical sample space }

The stochastic calculus involving   fBM is currently being
developed; see e.g. \cite{Nualart, Viens} and references therein.
This will lead to more advances   in the study of SDEs driven by
colored fBM noise:
\begin{equation}
dX_t   =  b(X_t )d t + \s(X_t )dB^{H}(t).
\end{equation}
Since the fBM  $B^{H}(t)$ is not Markov, the  solution process
$X_t$ is not Markov either. Thus the usual techniques from Markov
processes will not be applicable to the study of SDEs driven by
fBms. However, the random dynamical systems approach, as described
in \S \ref{RDS} above, looks promising \cite{MS, GMS}. The theory
of RDS, developed by Arnold and coworkers  \cite{Arnold},
describes the qualitative behavior of systems of  stochastic
differential equations in terms of stability, Lyapunov exponents,
invariant manifolds, and attractors.

See Theorem 2.3 in  \cite{MS, Nualart-flow}, the canonical sample
space is $\Om:=C_0(\R, \R^n)$, the set of continuous functions
that is zero at zero, but the probability measure $\mu_{fBM}$ is
generated by $B^H_t$ under the compact-open topology as defined in
Section \ref{whitenoise}. The Borel $\s-$field is $\F:=\B(C_0(\R,
\R^n))$.

We can introduce a   flow $\theta_t$ on this canonical sample
space $\Om$   defined by the   shifts
 \begin{equation*}
     \theta_t\omega(\cdot)=\omega(\cdot+t)-\omega(t),\quad
     t\in\mathbb{R},\qquad \omega\in\Omega=C_0(\mathbb{R},\mathbb{R}^d).
 \end{equation*}
 In this case the measure $\mu_{fBM}$ is invariant, i.e.
\begin{equation*}
\mu_{fBM}(\theta_{t}^{-1}(A))=\mu_{fBM}(A)
\end{equation*}
for all $A \in \mathcal{F}$.

\section{Dynamical systems driven by non-Gaussian noises}
\label{nongaussian}

In the last two sections, we considered dynamical systems driven
by Gaussian noises (white or colored), in terms of Brownian motion
or fractional Brownian motion. In this section, we discuss
  differential equation driven by non-Gaussian L\'evy
noises.

\subsection{L\'evy Motions  }

Gaussian processes like Brownian motion have been widely used to
model fluctuations in   engineering and science. For a particle in
Brownian motion, its sample paths are continuous in time almost
surely (i.e., no jumps), its mean square displacement increases
linearly in time (i.e., normal diffusion), and  its probability
density function decays exponentially in space (i.e., light tail
or exponential relaxation) \cite{Oksendal}.  But some complex
phenomena involve non-Gaussian fluctuations, with peculiar
properties such as anomalous diffusion (mean square displacement
is a nonlinear power law of time) \cite{BD90} and heavy tail
(non-exponential relaxation) \cite{Yon96}. For instance, it has
been argued that diffusion in a case of geophysical turbulence
\cite{Shlesinger} is anomalous. Loosely speaking, the diffusion
process consists of a series of ``pauses", when the particle is
trapped by a coherent structure, and ``flights" or ``jumps" or
other extreme events, when the particle moves in a jet flow.
Moreover, anomalous electrical transport properties have been
observed in some amorphous materials such as insulators,
semiconductors and polymers, where transient current is
asymptotically a power law function of time \cite{SSB91,
Herrchen}.  Finally, some paleoclimatic data \cite{Dit} indicates
heavy tail distributions and some DNA data \cite{Shlesinger} shows
long range power law decay for spatial correlation.

L\'evy motions are thought to be appropriate models for
non-Gaussian processes with jumps \cite{Sato-99}. Let us recall
that  a L\'evy motion $L(t)$, or $L_t$, is a non-Gaussian process
with independent and stationary increments, i.e., increments
$\Delta L (t, \Delta t)= L(t + \Delta t)-  L(t)$ are stationary
(therefore $\Delta L$ has no statistical dependence on $t$) and
independent for any non overlapping time lags $\Delta t$.
Moreover, its sample paths are only continuous in probability,
namely, $\PX (|L(t)-L(t_0)| \geq \delta) \to 0$ as $t\to t_0$ for
any positive $\delta$. With a suitable modification
\cite{Applebaum}, these paths may be taken as   c\`{a}dl\`{a}g,
i.e., paths are continuous on the right and have limits on the
left. This continuity is weaker than the usual continuity in time.

This generalizes the Brownian motion $B(t)$ or $B_t$, as $B(t)$
satisfies all these three conditions.  But \emph{Additionally},
(i) Almost every sample path of the Brownian motion     is
continuous in time in the usual sense and (ii) Brownian motion's
increments are Gaussian distributed.

\medskip

  Dynamical systems  driven by   non-Gaussian L\'evy noises
have attracted much attention  recently \cite{Applebaum, JW,
Schertzer}. Under certain conditions, the SDEs driven by L\'evy
motion   generate stochastic flows \cite{Kunita2004, Applebaum},
and also generate random dynamical systems (or cocycles)  in the
sense of Arnold \cite{Arnold}.  Recently,   exit time estimates
have been investigated by Imkeller and Pavlyukevich
\cite{ImkellerP-06, ImkellerP-08} , and Yang and Duan
\cite{YangDuan} for SDEs driven by L\'evy motion. This shows some
qualitatively different dynamical behaviors between SDEs driven by
Gaussian and non-Gaussian noises.

\medskip

  \emph{L\'evy motions} are named in honor of the
French probabilist Paul L\'evy, who first studied them in 1930s.
From a mathematical point of view, there
are so many reasons why they are so important \cite{Applebaum}, such as:\\
$\bullet$ There are many important examples, such as Brownian
motion, the Poisson process, stable processes, and subordinators.\\
$\bullet$ They are generalizations of random walks to continuous
time.\\
$\bullet$ They are the simplest class of processes whose paths
consist of continuous motion interspersed with jump
discontinuities
of random size appearing at random times.\\

\noindent\textbf{Definition} \cite{Applebaum} A stochastic process
$L=(L(t),t\geq 0)$ defined on a probability space $(\Omega,
\mathcal {F},P)$ is a L\'evy motion if:\\
(\textbf{L1}) $L(0)=0$ a.s.; \\
(\textbf{L2}) $L$ has independent and stationary increments; \\
(\textbf{L3}) $L$ is \emph{stochastically continuous}, i.e. for
all $\epsilon>0$ and for all $s$
\begin{equation*}
\lim_{t\rightarrow s}P(|X(t)-X(s)|>\epsilon)=0
\end{equation*}
With a suitable modification \cite{Applebaum}, these paths may be
taken as   c\`{a}dl\`{a}g, i.e., paths are continuous on the right
and have limits on the left. This continuity is weaker than the
usual continuity in time.

Figure 3 is a sample path for a L\'evy motion.
\begin{figure} \label{levypath}
\includegraphics{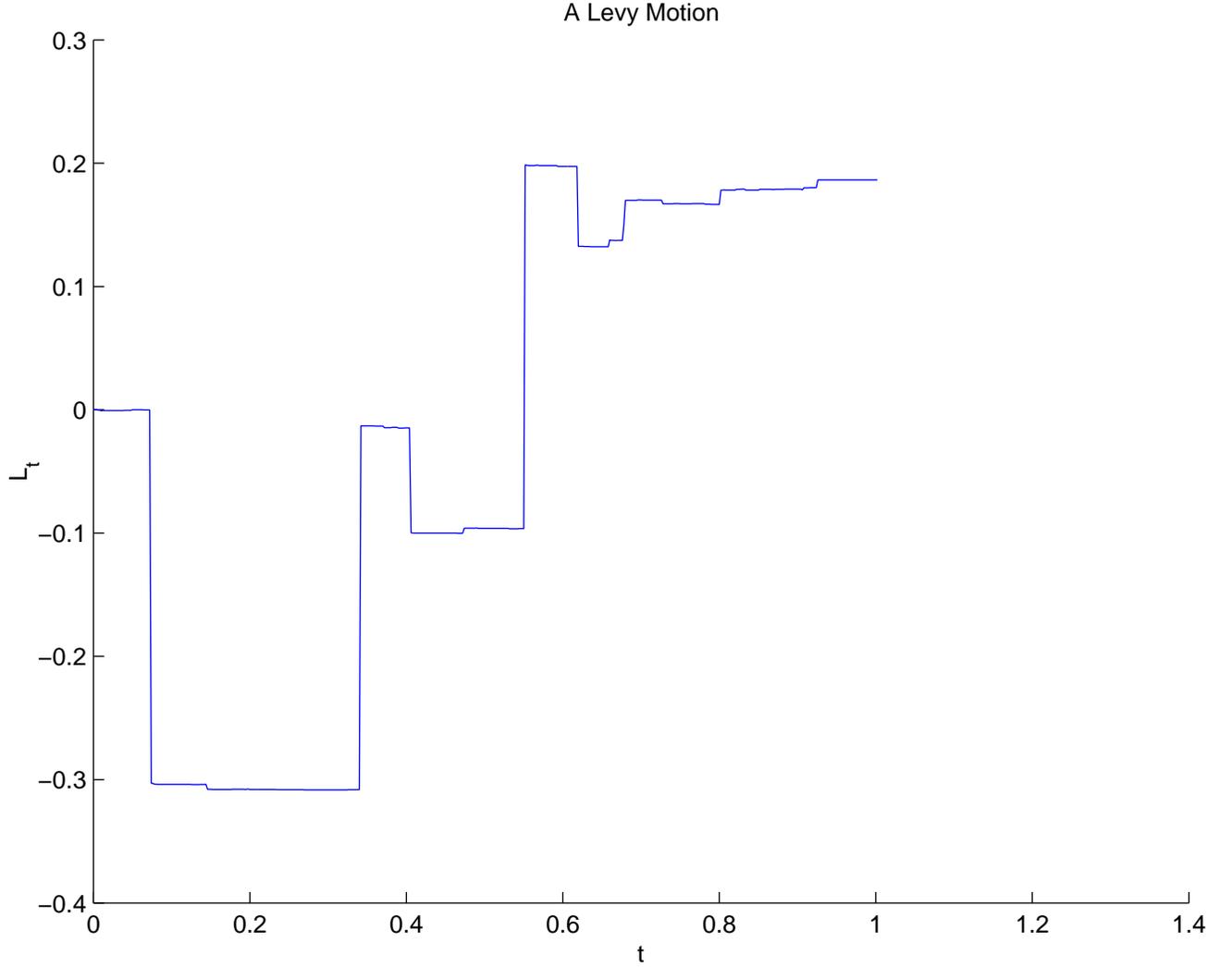}
\caption{A sample path for a L\'evy motion}
\end{figure}

One way to understand the structure of the L\'evy motions is to
employ Fourier analysis. It may be shown that   each $L(t)$ is
\emph{infinitely divisible}. The infinitely divisible random
variables are characterized completely through their
characteristic functions by a beautiful formula,   established by
Paul L\'evy and A. Ya. Khintchine in the 1930s. The
related definitions and theorems are as follows.\\
\\
\noindent\textbf{Definition} \cite{Applebaum} The
\emph{characteristic function}  of a stochastic process $X(t)$
taking values in $\mathbb{R}^{d}$ is the mapping $\Phi_{t}:
\mathbb{R}^{d}\rightarrow\mathbb{C}$ defined by
\begin{equation*}
\Phi_{t}(u)=\mathbb{E}(e^{iu\cdot
X(t)}):=\int_{\mathbb{R}^{d}}e^{iu\cdot y}p_{t}(dy)
\end{equation*}
where $p_{t}$ is the distribution of $X(t)$.\\
\\
\noindent\textbf{Definition} \cite{Applebaum} $X(t)$ is
\emph{infinitely divisible} if for each $n\in\mathbb{N}$, there
exists a probability measure $p_{t,n}$ on $\mathbb{R}^{d}$ with
characteristic function $\Phi_{t,n}$ such that
$\Phi_{t}(u)=(\Phi_{t,n}(u))^{n}$, for each $u\in \mathbb{R}^{d}$.\\

\begin{theorem}[The L\'evy-Khintchine Formula \cite{Applebaum}]
If $L=(L(t),t\geq0)$ is a L\'evy motion, then
\begin{equation*}
\Phi_{t}(u)=e^{t\eta(u)}\qquad for\ each\ t\geq0,
 u \in  \mathbb{R}^{d}
\end{equation*}
where
\begin{equation}\label{eq4}
\eta(u)=ib\cdot u-\frac{1}{2}u\cdot
au+\int_{\mathbb{R}^{d}-\{0\}}[e^{iu\cdot y}-1-iu\cdot y \;
\mathbf{I}_{|y|<1}(y)]\nu(dy)
\end{equation}
for some $b\in\mathbb{R}^{d}$, a non-negative definite symmetric
$d\times d$ matrix $a$ and a Borel measure $\nu$, called L\'evy
jump measure, on $\mathbb{R}^{d}-\{0\}$ for which
$\int_{\mathbb{R}^{d}-\{0\}}(|y|^{2}\wedge1)\nu(dy)<\infty$.
Here $\mathbf{I}_S$ is the indicator function of the set $S$.\\
Conversely, given a mapping of the form \eqref{eq4}, we can always
construct a L\'evy motion for which $\Phi_{t}(u)=e^{t\eta(u)}$ and
call it a L\'evy motion with the characteristics or generating
triple $(a,b,\nu)$.
\end{theorem}

\textbf{Remark}: The so-defined Borel measure $\nu$ in \eqref{eq4}
is called the \emph{L\'evy jump measure}, which should not be
confused with the
\emph{probability  measure} induced by the L\'evy motion, to be defined below.\\
\\
The different terms which appear in the L\'evy-Khintchine formula
have a probabilistic significance emphasized in \cite{RevuzYor}.
Every L\'evy motion  is obtained as a sum of independent processes
with three types of characteristics $(0,b,0)$, $(a,0,0)$ and
$(0,0,\nu)$. Thus, the L\'evy measure accounts for the jumps of
$L$ and the knowledge of $\nu$ permits to give a probabilistic
construction of $L$; see \cite{RevuzYor}.

\subsection{Canonical sample space}

Consider a SDE  driven by non-Gaussian Levy noise
\begin{equation}
dX_t   =  b(X_t )d t + \s(X_t )dL(t).
\end{equation}

The canonical space  has to be  enlarged to include all the
\emph{cadlag} functions, i.e. functions that are right-continuous
and have left limits,  defined on $\mathbb{R}$ and taking values
in $\mathbb{R}^d$. This space is denoted as $D(\R, \R^n)$.


We adopt the same point of view as in Section \ref{whitenoise}
that a stochastic process is also a
\emph{random variable}, i.e.\\
\begin{equation*}
L_{t}(\omega): \omega\rightarrow D(\R, \R^n)
\qquad\omega\rightarrow L(t), t\in \R.
\end{equation*}

\begin{remark}  Here the compact-open metric defined in \eqref{eq2}
cannot make $D(\R^+, \R^n)$   separable \cite{Protter}. For
example, if
$f_{\alpha}(t)=1_{[\alpha,\infty)}(t)$,$f_{\beta}(t)=1_{[\beta,\infty)}(t)$,
then $\rho(f_{\alpha},f_{\beta})=1/2$ for all $\alpha, \beta$ with
$0\leq\alpha<\beta\leq1$, and since there are uncountably many
such $\alpha,\beta$, the space is not separable.
\end{remark}

However, it can be made   complete and separable when endowed with
the \emph{Skorohod metric}   \cite{Billingsley}. With this special
metric, we call $D(\mathbb{R},\mathbb{R}^d)$ a Skorohod space. The
Skorohod metric on $D(\mathbb{R},\mathbb{R}^d)$ is defined as
\begin{equation*}
d(x,y):=\sum_{m=1}^{\infty}\frac{1}{2^{m}}(1\wedge
d_{m}^{\circ}(x^{m},y^{m}))\qquad for\ all \ x, y\in D
\end{equation*}
where $x^{m}(t):=g_{m}(t)x(t)$, $y^{m}(t):=g_{m}(t)y(t)$ with
\begin{equation*}
g_{m}(t):= \left\{
\begin{array}{rl}
1,& \text{if } |t|\leq m-1
\\[2ex]
m-t,& \text{if } m-1\leq |t| \leq m,
\\[2ex]
0,& \text{if } |t|\geq m
\end{array} \right.
\end{equation*}
and
\begin{equation*}
d_{m}^{\circ}(x,y):=\inf_{\lambda\in\Lambda}\left\{\sup_{-m\leq
s<t\leq m}\left|\log\frac{\lambda (t)-\lambda (s)}{t-s}\right|\vee
\sup_{-m\leq t\leq m}|x(t)-y(\lambda (t))|\right\},
\end{equation*}
where $\Lambda$ denotes the set of strictly increasing, continuous
functions from $\R$ to itself.\\

The Borel $\sigma-$field under this topology is denoted as
$\mathcal{S}$.
For studying the \emph{weak convergence}
and \emph{tightness} in $D$, the same approach adopted in $C$ can
be applied except that the fact the natural projections are not
continuous need to be noticed \cite{Billingsley}.\\

\noindent\textbf{Definition} The probability measure, $\mu_L$, in
$(D(\R, \R^n), \mathcal{S})$ that makes every element in $D(\R,
\R^n)$ a sample L\'evy path is called the L\'evy probability
measure. Note that this measure is not to be confused with the
L\'evy jump measure $\nu$ mentioned above.


We can also introduce a flow $\theta=(\theta_{t},t\in \R) $ on
this canonical sample space $\Omega$   by the shifts
\begin{equation}\label{shift}
(\theta_{t}\omega)(s):=\omega(t+s)-\omega(t).
\end{equation}
The (L\'evy) probability measure $\mu_L$ is  invariant under this
flow, i.e.
\begin{equation*}
\mu_{L}(\theta_{t}^{-1}(A))=\mu_{L}(A)
\end{equation*}
for all $A\in\mathcal {F}$; see \cite{Applebaum} (Page 325) and
\cite{Liu}. This flow is an ergodic dynamical system \cite{Arnold}
with respect to the above probability measure $\mu_L$.


\end{document}